\documentclass[12pt]{amsart}

\usepackage{amssymb}
\usepackage{latexsym}
\usepackage{verbatim}
\usepackage{fullpage}
\usepackage{pslatex}

\input {cyracc.def}
 \font\tencyr=wncyr10 
\font\tencyi=wncyi10 
\font\tencysc=wncysc10 
\def\rus{\tencyr\cyracc}
\def\rusi{\tencyi\cyracc}
\def\rusc{\tencysc\cyracc}

\renewenvironment{proof}
{\noindent {\sl Proof.}\quad }{\hfill
$\square$ \vskip1.1ex\noindent }

\newenvironment{proof*}
{\noindent {\sl Proof.}\quad }{\hfill
$\square$}

\renewcommand{\theequation}{\thesection .\arabic{equation}}
\renewcommand{\thesubsubsection}{\theequation .\arabic{subsubsection}}

\catcode`\@=\active
\catcode`\@=11
\def\@eqnnum{\hbox to
.01pt{}\rlap{\hskip-\displaywidth(\mathbf{\theequation})}}
\catcode`\@=12

\newenvironment{s}[1]
{ \vskip1.2ex \refstepcounter{equation}
\noindent {\bf \theequation\enspace #1.} \begin{sl}}{\end{sl}
\vskip1.1ex\noindent }

\newenvironment{rem}[1]
{ \vskip1.2ex \refstepcounter{equation}
\noindent {\bf \theequation\enspace {#1}.} }{ \vskip1.1ex\noindent }

\newenvironment{subs}[1]
{\vskip1.2ex \refstepcounter{equation}
\noindent {\bf (\theequation)\quad #1.} }{\quad}

\newcommand {\be}{{\frak b}}

\newcommand {\g}{{\frak g}}

\newcommand {\el}{{\frak l}}

\newcommand {\te}{{\frak t}}

\newcommand {\gC}{{\goth C}}
\newcommand {\gH}{{\goth H}}

\newcommand {\gR}{{\goth R}}

\newcommand {\ap}{\alpha}

\newcommand {\cO}{{\mathcal O}}

\newcommand {\VV}{{\Bbb V}}

\newcommand {\Lie}{{\mathrm{Lie\,}}}

\newcommand {\rk}{{\mathrm{rk\,}}}

\newcommand {\tri}{{\frak sl}_2}
\newcommand {\GR}[2]{{\textrm{{\bf #1}}}_{#2}}

\newcommand {\ov}{\overline}

\newcommand {\pv}{{\mathcal P}^*(\VV)}

\newcommand {\beq}{\begin{equation}}
\newcommand {\eeq}{\end{equation}}

\newfam\Bbbfam\newfam\eufam\newfam\eusfam%
\font\Bbbfont=msbm10 scaled 1200%
\font\olala=msam10 scaled 1200%
\font\frak=eufm10 scaled 1400%
\font\Bbbsmallfont=msbm8%
\textfont\Bbbfam=\Bbbfont\scriptfont\Bbbfam=\Bbbsmallfont
\font\euzw=eufm10 scaled 1200%
\font\euac=eufm7 scaled 1200%
\font\euacc=eufm7 scaled 1000%
\textfont\eufam=\euzw\scriptfont\eufam=\euac%
\scriptscriptfont\eufam=\euacc%
\font\euszw=eusm10 scaled 1200%
\font\eusac=eusm7 scaled 1200%
\font\eusacc=eusm7 scaled 1000%
\textfont\eusfam=\euszw\scriptfont\eusfam=\eusac%
\scriptscriptfont\eusfam=\eusacc%
\def\frak{\fam\eufam}%
\def\goth{\fam\eusfam}%
\def\Bbb{\fam\Bbbfam}%

\def\varnothing{\hbox {\Bbbfont\char'077}}
\def\square{\hbox {\olala\char"03}}

\begin{document}
\setlength{\parskip}{2pt plus 4pt minus 0pt}
\hfill {\scriptsize June 6, 2003}
\vskip1ex
\vskip1ex

\title[Regions and nilpotent orbits]{Regions in the dominant
chamber and nilpotent orbits}
\author{\sc Dmitri I. Panyushev}
\thanks{This research was supported in part by R.F.B.I. Grants no.
01--01--00756 and  02--01--01041}
\maketitle
\begin{center}
{\footnotesize
{\it Independent University of Moscow,
Bol'shoi Vlasevskii per. 11 \\
121002 Moscow, \quad Russia \\ e-mail}: {\tt panyush@mccme.ru }\\
}
\end{center}

\noindent
Let $G$ be a complex semisimple algebraic group with Lie algebra $\g$.
The goal of this note is to show that combining some ideas of
\cite{VP} and \cite{gs} quickly yields a geometric
description of the characteristic of a nilpotent $G$-orbit in
an arbitrary (finite-dimensional) rational $G$-module.
\\
Fix a Borel subalgebra $\be\subset\g$ and a Cartan subalgebra $\te$ in it.
For each $\te$-weight $\mu$ of a $G$-module $\VV$, consider the affine
hyperplane $\gH_{\mu,2}=\{x \mid (x,\mu)=2\}\subset \te_{\Bbb R}$.
These hyperplanes cut the dominant chamber in finitely many
regions, and to any region
$R$ one may attach a $\be$-stable subspace of $\VV$ by the following rule:
\[
     \VV_R=\bigoplus_{\mu\in I_R}\VV^\mu \ ,
\]
where $I_R$ is the set of weights of $\VV $ such that $(x,\mu)>2$
for some (equivalently, any) $x\in R$. Given a nilpotent $G$-orbit
$\cO\subset\VV$, consider the closure of the union of all regions $R$ such
that $\cO\cap \VV_R\ne\varnothing$. Let's call this set $\gC_\cO$.
Our first observation is that $\gC_\cO$ contains a unique element
of minimal length, and this element is just the dominant characteristic
of $\cO$ in the sense of \cite[5.5]{VP}. Next, we show that if
the representation $G\to GL(\VV)$ is associated with either a periodic or
a ${\Bbb Z}$-grading of a reductive algerbaic Lie algebra, then
the condition ``$\cO\cap \VV_R\ne\varnothing$'' can be replaced with
``$\cO\cap \VV_R$ is dense in $\VV_R$''. This new condition determines
a smaller set $\tilde\gC_\cO\subset \gC_\cO$, but
these two sets still have the same element of minimal length.
This provides another proof and also a generalization of the main result
of \cite{gs}. It is worth noting that the representations associated
with ${\Bbb Z}_m$-gradings are {\it visible\/}, i.e., contain finitely many nilpotent orbits, and in this case different orbits have different characteristics.
\\
We also give an example showing that, for an arbitrary visible $G$-module
$\VV$, it may happen that different orbits have the same characteristic
and that for some orbits $\cO$ there are no subspaces of the form $\VV_R$
such that $\cO\cap \VV_R$ is dense in $\VV_R$.



\noindent
\begin{subs}{Main notation}
\end{subs}
$\Delta$ is the root system of $(\g,\te)$ and
$W$ is the Weyl group of $(\te,\Delta)$.

$\Delta^+$  is the set of positive
roots and
$\Pi=\{\ap_1,\dots,\ap_p\}$ is the set of simple roots in $\Delta^+$.
 \\
We define $\te_{\Bbb R}$ to be set of all elements of $\te$ having real
eigenvalues in any $G$-module (a Cartan subalgebra of a split real form
of $\g$). Denote by $(\ ,\ )$ a $W$-invariant inner product on
$\te_{\Bbb R}$. Using $(\ ,\ )$, we identify $\te_{\Bbb R}$ and
$\te_{\Bbb R}^*$. So that, one may think that
$\te_{\Bbb R}=\oplus_{i=1}^p{\Bbb R}\ap_i$.


${\gC}=\{x\in V\mid (x,\ap)>0 \ \ \forall \ap\in\Pi\}$
\ is the (open) fundamental Weyl chamber.




\section{The characteristic of a nilpotent orbit}
\label{vinberg}
\setcounter{equation}{0}

\noindent
In this section we recall some results published
in the survey article~\cite[\S 5]{VP}. Unfortunately, that simple approach
to questions of stability, optimal one-parameter subgroups, and a
stratification of the null-cone remained
largely unnoticed by the experts.
\\[.5ex]
Let $\VV$ be a $G$-module. Write $\VV^\mu$ for the $\mu$-weight space
of $\VV$. Here $\mu$ is regarded as element of $\te_{\Bbb R}$. Hence
$\mu(x)=(\mu,x)$ for any $x\in \te_{\Bbb R}$.
\\[.5ex]
Suppose $h\in\te_{\Bbb R}$, i.e., $h$ is a rational semisimple element.
For a $G$-module $\VV$ and $c\in {\Bbb Q}$,
we set
\[
\VV_h\langle c\rangle =\{ v\in\VV \mid h{\cdot}v=cv\},\quad
\VV_h\langle {\ge} c\rangle =\oplus_{k\ge c} \VV_h\langle k\rangle ,
\text{ \ and \ }
\VV_h\langle {>} c\rangle =\oplus_{k > c} \VV_h\langle k\rangle  \ .
\]
For instance, $\g_h\langle 0\rangle $ is the centralizer of $h$ in $\g$
(a Levi subalgebra of $\g$), $\g_h\langle{\ge} 0\rangle$ is a
parabolic subalgebra of $\g$, and
$\g_h\langle {>} 0\rangle $ is the nilpotent radical of
$\g_h\langle{\ge} 0\rangle$.
Clearly,
\[
   \VV_h\langle c\rangle =
\bigoplus_{\mu:\, \mu(h)=c}\VV^\mu\quad \text{ and }\quad
  \g_h\langle a\rangle {\cdot}\VV_h\langle c\rangle
  \subset \VV_h\langle a+c\rangle  \ .
\]
Recall that an element $v\in\VV$, or the orbit $G{\cdot}v$,
is called {\it nilpotent},
if $\ov{G{\cdot}v}\ni 0$. It is easy to verify that, for any
$h\in\te_{\Bbb R}$, the subspace
$\VV_h\langle {>}0\rangle $ consists of nilpotent elements.
Conversely, the Hilbert-Mumford criterion asserts that any nilpotent
$G$-orbit in $\VV$ meets a subspace of this form for a suitable $h$.
\begin{rem}{Definition}   \label{def-char}
The {\it characteristic\/} of a nilpotent orbit $\cO$ is a shortest
element $h\in\te_{\Bbb R}$ such that $\cO\cap \VV_h\langle{\ge} 2\rangle\ne
\varnothing$.
\end{rem}%
{\sf Remark.}
In principle, one may choose an arbitrary normalization "$(\ge c)$"
in the definition. The choice $c=2$ is explained by the fact that
for $\VV=\g$ this leads to the usual (Dynkin) characteristic of a nilpotent
element.
\\
It was shown in \cite[5.5]{VP} that each nilpotent orbit has a
characteristic. Moreover, if $h_1,h_2\in \te_{\Bbb R}$ are two
characteristics of $\cO$, then they are $W$-conjugate.
Thus, to any nilpotent orbit $\cO\subset\VV$
one may attach uniquely the {\it dominant
characteristic\/}, which is denoted by $h_\cO$.
\\
If we are given an $h\in\te_{\Bbb R}$ and $u\in \VV_h\langle {\ge} 2\rangle$,
then it is
helpful to have a criterion to decide whether $h$ is a characteristic
of $G{\cdot}u$. The following result, attributed in \cite[Theorem\,5.4]{VP}
to F.\,Kirwan and L.\,Ness, solves the problem.
Let $Z_G(h)$ denote the centralizer of $h$ in $G$ and $\tilde Z_G(h)
\subset Z_G(h)$ the {\it reduced centralizer\/}.
That is, the Lie algebra of $\tilde Z_G(h)$
is the orthogonal complement to $h$ in $\g_h\langle 0\rangle=\Lie Z_G(h)$.
Clearly, $\VV_h\langle c\rangle$ is a $Z_G(h)$-module for any $c$.

\begin{s}{Theorem}  \label{crit}
Under the previous notation, $h$ is a characteristic of $G{\cdot}u$
if and only if the projection of $u$ to $\VV_h\langle 2\rangle$ is not a nilpotent element
with respect to the action of $\tilde Z_G(h)$.
\end{s}%
%


\section{Regions and characteristics}
\label{R&C}
\setcounter{equation}{0}

\noindent
Let $\VV$ be a $G$-module. Write
${\mathcal P}^*(\VV)$ for the set of nonzero weights of
$\VV$ with respect to $\te$. For any
$\mu\in{\mathcal P}^*(\VV)$, consider the affine
hyperplane $\gH_{\mu,2}=\{x\in \te_{\Bbb R}\mid (x,\mu)=2\}$.
The number ``2" is determined by the normalization in
Definition~\ref{def-char}.
We will be interested in the hyperplanes meeting the dominant Weyl chamber.
It is easily seen that the following is true.
\begin{s}{Lemma} \label{hyper}
We have \\
\indent $\gH_{\mu,2}\cap \gC\ne \varnothing \ \Longleftrightarrow$ \
$\mu$ has a positive
coefficient in the expression
$\mu=\displaystyle\sum_{i=1}^p a_i\ap_i$ $(a_i\in {\Bbb Q})$.
\end{s}%
The set of all such hyperplanes cuts $\gC$ in regions.  That is,
a {\it region\/} (associated with $\VV$) is a connected component
of $\gC\setminus \displaystyle\bigcup_\mu \gH_{\mu,2}$. The set of all
regions is denoted by $\gR=\gR(\VV)$.
Clearly, the closure of each region is a convex polytope.
Given a region $R$, consider all hyperplanes separating $R$ from the origin,
and the corresponding weights in ${\mathcal P}^*(\VV)$. This set
of weights is denoted by $I_R$. More precisely, if $x\in R$, then
\[
  I_R=\{\mu\in \pv \mid (x,\mu)>2\} \ .
\]
\vskip-1.5ex
\begin{s}{Lemma} Let $R\in \gR$.
\begin{itemize}
\item[\sf (i)] \ 
if $\mu\in I_R$, $\gamma\in\Delta^+$, and $\mu+\gamma\in \Delta^+$,
then $\mu+\gamma\in I_R$.
\item[\sf (ii)] \ The subspace $\VV_R:=\oplus_{\mu\in I_R} \VV^\mu
\subset \VV$ is $\be$-stable.
\item[\sf (iii)] \ Each $G$-orbit meeting $V_R$ is nilpotent.
\end{itemize}
\end{s}\begin{proof}
(i) - obvious; (ii) follows from (i); (iii) we have $\displaystyle
\lim_{t\to -\infty}
\exp(tx){\cdot}u=0$ for any $x\in R$ and $u\in V_R$.
\end{proof}%
Suppose $\cO\subset\VV$ is a nilpotent $G$-orbit.
One may attach to $\cO$ a collection of regions,  as follows.
Set
\begin{equation} \label{sootv}
  M_\cO= \{R\in\gR \mid \cO\subset G{\cdot}\VV_R\}=
        \{R\in\gR \mid \cO\cap\VV_R\ne\varnothing \} \ ,
\end{equation}
and
\[
  \gC_\cO=\bigcup_{R\in M_\cO} \bar R \subset \bar\gC \ .
\]
Thus, $\gC_\cO$ is a closed subset of $\bar\gC$ determined by
$\cO$. Let $h' \in \gC_\cO$ be an element of minimal
length.

\begin{s}{Proposition}
$h'$ is a unique element of minimal length in $\gC_\cO$,
and $h'=h_\cO$.
\end{s}\begin{proof}
By the very construction, $h'$ has the property
that $\cO\cap \VV_{h'}\langle {\ge 2}\rangle \ne \varnothing $ and it is a
shortest dominant element with this property. It then follows from results
described in Section~\ref{vinberg} that
$h'=h_\cO$ and $\gC_\cO$ contains a unique element
of minimal length.
\end{proof}%
The above construction is inspired by \cite{gs},
where the case $\VV=\g$ is considered.
However, the condition \ref{sootv} was slightly different there. Namely,
the set of regions attached to $\cO$ was determined by the condition
that $\cO\cap \g_R$ be dense in
$\g_R$. But this stronger condition cannot lead in general
to satisfactory results,
unless $\VV$ is a {\it visible\/} $G$-module.
For, the number of subspaces of the form
$\VV_R$ is finite and therefore the set of such regions would be empty for
infinitely many nilpotent orbits, if $\VV$ is not visible.
Moreover, even if $\VV$ is visible, it may happen
that, for a given nilpotent orbit, there is no subspace
$\VV_R$ $(R\in\gR)$ such that $\cO\cap \VV_R$ is dense in $\VV_R$
(see example below).

However, one may formally set
\begin{equation} \label{visible}
 \widetilde M_\cO= \{R\in\gR \mid \cO \ \text{ is dense in }\  G{\cdot}\VV_R\}=
        \{R\in\gR \mid \cO\cap\VV_R \ \text{ is dense in }\ \VV_R \} \ ,
\end{equation}
and
\[
  \widetilde\gC_\cO=\bigcup_{R\in \tilde M_\cO} \bar R \subset \bar\gC \ .
\]
Clearly, $\widetilde M_\cO\subset M_\cO$ and $\widetilde\gC_\cO\subset
\gC_\cO$. We also define $\tilde h_\cO$ to be an element of
$\widetilde\gC_\cO$ of minimal length
(if $\widetilde\gC_\cO\ne\varnothing$ !)

\begin{rem}{Example}
Here we give an example of a visible module such that
(i) $\widetilde M_\cO=\varnothing$ for some $\cO$, and (ii)
$h_{\cO_1}=h_{\cO_2}$ for different nilpotent orbits.

Let $G=SL(V_1)\times SL(V_2)$, $\dim V_1=\dim V_2=2$, and
$\VV=(V_1\otimes V_2) \oplus V_1$.
\\
Identifying $\VV$ with the space of 2 by 3 matrices, we write
$v=\left( \begin{array}{ccc}
m & n & x \\ p & q & y
\end{array}\right)$ for a generic element in $\VV$.
Here $V_1\otimes V_2=\left\{\left( \begin{array}{ccc}
m & n & 0 \\ p & q & 0
\end{array}\right)\right\}$.
There are five nilpotent orbits in $V$,
and representatives of the non-trivial orbits are:

$\cO_2$ :  \ $v_2=\left( \begin{array}{ccc}
0 & 0 & 1 \\ 0 & 0 & 0  \end{array}\right)$,
\qquad
$\cO_3$ : \  $v_3=\left( \begin{array}{ccc}
1 & 0 & 0 \\ 0 & 0 & 0  \end{array}\right)$,

$\cO_4$ :  \ $v_4=\left( \begin{array}{ccc}
1 & 0 & 1 \\ 0 & 0 & 0  \end{array}\right)$,
\qquad
$\cO_5$ : \  $v_5=\left( \begin{array}{ccc}
1 & 0 & 0 \\ 0 & 0 & 1  \end{array}\right)$.

\noindent
One has $\dim\cO_i=i$ for $2\le i\le 5$.
Since $\rk G=2$, we have $\te_{\Bbb R}$ is ${\Bbb R}^2$
and it is not hard to depict all the regions
and determine the characteristics.
\begin{figure}[htb]
 \label{hasse}
\setlength{\unitlength}{0.02in}
\centerline
{
\begin{picture}(50,100)(15,-15)
\put(0,0){\vector(0,1){80}} \put(4,75){$x_2$}
\put(0,0){\vector(1,0){100}} \put(95,5){$x_1$}
\put(4,4){I}
\put(8,18){II}    \put(20,40){\circle{4}}
\put(5,40){III}    \put(20,0){\circle{4}}
\put(33,4){IV}    \put(10,10){\circle{4}}
\put(25,25){V}
\put(26,62){VI}
\thicklines
\put(20,-10){\line(0,1){80}}
\put(-10,10){\line(1,1){55}}
\put(10,-10){\line(1,1){60}}
\put(30,-10){\line(-1,1){50}}
\end{picture}
}
\end{figure}
The dominant chamber is the positive quadrant. There are four lines
of the form $\gH_{\mu,2}$ meeting
the positive quadrant, which correspond to the roots
$\ap_1,\,\ap_1+\ap_2$,\,$\ap_1-\ap_2,\,-\ap_1+\ap_2$. Here
$\ap_i$ is the simple root of $SL(V_i)$. Hence
one gets six regions marked by Roman numbers. We have
$\tilde M_{\cO_2}=\varnothing$,\
$\tilde M_{\cO_3}=\{\text{II,\,III}\}$,\
$\tilde M_{\cO_4}=\{\text{IV,V}\}$,\
$\tilde M_{\cO_5}=\{\text{VI}\}$.
Since  $\ov{\cO_4}\supset \cO_2$ and  $\ov{\cO_3}\not\supset\cO_2$, we
conclude that $h_3=(1,1)$, $h_2=h_4=(2,0)$, and $h_5=(2,4)$.
The elements $h_i$ are circled in the figure.
\end{rem}%
Therefore one should not expect that $\tilde h_\cO$ is always
defined and that different orbits have different characteristics.
At the rest of the section,
we give a sufficient condition for this to happen.

Let $\el$ be a reductive algebraic Lie algebra.
Consider a ${\Bbb Z}_m$-grading of $\el$, where $m\in \Bbb N$ or
$m=\infty$. That is, we have
$\el=\underset{i\in{\Bbb Z}_m}{\bigoplus}\el_i$ if $m$ is finite,
and $\el=\underset{i\in{\Bbb Z}}{\bigoplus}\el_i$ is a
$\Bbb Z$-grading in the second case. Here $\el_0$ is reductive and
each $\el_i$ is an $\el_0$-module.
Let $G$ be a connected (reductive) group with Lie algebra $\el_0$,
and set $\VV=\el_1$. Then we shall say that the representation
$G\to GL(\VV)$ is associated with a ${\Bbb Z}_m$-grading (of $\el$).
By a famous result of
Vinberg \cite[\S\,2]{vi76}, $\VV$ is a visible $G$-module in
this situation.

\begin{s}{Theorem}
Suppose the representation
$G\to GL(\VV)$ is associated with a ${\Bbb Z}_m$-grading.
Then different nilpotent $G$-orbits in $\VV$ have different characteristics
and for any nilpotent $G$-orbit $\cO\subset\VV$ we have
$\tilde h_\cO=h_\cO$.
\end{s}\begin{proof}
1. Let $e\in\VV=\el_1$ be a nilpotent element, and $\cO=G{\cdot}e$.
By a generalization of the Morozov-Jacobson theorem
\cite[Theorem\,1(1)]{vi79}, there is an
$\tri$-triple $(e,h,f)$ such that $h\in \el_0=\g$ and $f\in\el_{-1}$.
The rational semisimple element $h$ determines a ${\Bbb Z}$-grading of
$\el$,
and we have $e\in \VV_h\langle 2\rangle \subset\el_h\langle 2\rangle$.
It is well known that, in the Lie algebra $\el$, we have
$\tilde Z_L(h){\cdot}e$ is closed in $\el_h\langle 2\rangle$.
Hence, by the Richardson-Vinberg lemma \cite[\S\,2]{vi76},
$\tilde Z_G(h){\cdot}e$ is closed in $\VV_h\langle 2\rangle $.
Without loss of generality, one may assume that $h$ is a dominant element
in $\te_{\Bbb R}$.
Then, in view of Theorem~\ref{crit},
$h=h_\cO$ is the dominant characteristic of $\cO$. Next,
$L{\cdot}e\cap \el_h\langle{\ge} 2\rangle$ is dense in
$\el_h\langle {\ge }2\rangle$ and hence
$\cO\cap \VV_h\langle{\ge} 2\rangle$ is dense in
$\VV_h\langle {\ge }2\rangle$.
Thus, $h=\tilde h_\cO$.

2. That different nilpotent orbits have different characteristics
stems from \cite[Theorem\,1(4)]{vi79}.
\end{proof}%
This result applies, in particular, to the adjoint representations
($m=1$), where one obtains another proof for the main result in
\cite{gs}.
Another interesting case is that of the little adjoint $G$-module,
if $\g$ is a simple Lie algebra having roots of different lengths.
Let $\theta_s$ be the short dominant root.
Then the simple $G$-module with highest weight $\theta_s$ is called
the {\it little adjoint\/}.
It is denoted by $\g_{la}$.
The set ${\mathcal P}^*(\g_{la})$ is $\Delta_s$, the set of short roots.
This module is always associated with a ${\Bbb Z}_m$-grading
($m=2$ for $\GR{B}{p}$, $\GR{C}{p}$, $\GR{F}{4}$; $m=3$ for $\GR{G}{2}$).
Therefore the set of regions $\gR(\g_{la})$ allows us to determine the
characteristics of the nilpotent $G$-orbits in $\g_{la}$.
The arrangement of hyperplanes $\gH_{\mu,2}$ ($\mu\in\Delta^+_s$) inside
of $\gC$ was
studied in \cite{short}, where it was shown that there is
a bijection between the set of regions $\gR(\g_{la})$ and the set of
all $\be$-stable subspaces of $\g_{la}$ without semisimple elements.
We also give in \cite{short} an explicit formula for the number
$\#\gR(\g_{la})$.




\end{document}